\documentclass[11pt,leqno]{article}

\usepackage{amsthm,amsfonts,amssymb,amsmath,color}
\usepackage{amsbsy}

\numberwithin{equation}{section}

\renewcommand\d{\partial}

\renewcommand\b{\beta}
\renewcommand\o{\omega}

\newcommand\R{\mathbb R}\newcommand\N{\mathbb N}

\def\g{\gamma}
\def\de{\mathfrak{z}}

\def\O{\Omega}
\def\th{\theta}

\def\l{\lambda}

\def\epsilon{\varepsilon}
\def\e{\varepsilon}

\newcommand\br{\begin{rem}}
\newcommand\er{\end{rem}}
\newcommand\bp{\begin{pmatrix}}
\newcommand\ep{\end{pmatrix}}
\newcommand\be{\begin{equation}}
\newcommand\ee{\end{equation}}
\newcommand\ba{\begin{equation}\begin{aligned}}
\newcommand\ea{\end{aligned}\end{equation}}

\newcommand\nn{\nonumber}

\usepackage{layout}

\newcommand{\TT}{{\mathbb T}}

\newcommand{\II}{{\mathbb I}}

\newcommand{\SSS}{{\mathbb S}}

\newcommand{\ff}{{\mathbf f}}

\newcommand{\vvarphi}{{\boldsymbol \varphi}}

\newcommand{\Ov}[1]{\overline{#1}}

\newcommand{\DC}{C^\infty_c}

\newcommand{\vr}{\varrho}

\newcommand{\vu}{\vc{u}}

\newcommand{\vc}[1]{{\bf #1}}
\newcommand{\Div}{{\rm div}}
\newcommand{\Grad}{\nabla_x}

\newcommand{\dx}{{\rm d} {x}}
\newcommand{\dt}{{\rm d} t }

\newcommand{\intO}[1]{\int_{\O} #1 \, \dx}

\newcommand{\dive}{{\rm div\,}}

\newtheorem{definition}{Definition}[section]
\newtheorem{theorem}[definition]{Theorem}
\newtheorem{proposition}[definition]{Proposition}

\newtheorem{remarks}[definition]{Remarks}

\begin{document}

\title{Global existence of large data weak solutions for a simplified compressible Oldroyd--B model without stress diffusion}

\author{Yong Lu \footnote{Department of Mathematics, Nanjing University, 210093 Nanjing, China. Email: {\tt luyong@nju.edu.cn}.} \and  Milan Pokorn\'y\footnote{Charles University, Faculty of Mathematics and Physics, Sokolovsk\'a 83, Prague, Czech Republic. Email: {\tt pokorny@karlin.mff.cuni.cz}.}  }

\date{}

\maketitle

\begin{abstract}

We start with the compressible Oldroyd--B model derived in \cite{Barrett-Lu-Suli} ({\em J. W. Barrett, Y. Lu, E. S\"uli.  Existence of large-data finite-energy global weak solutions to a compressible Oldroyd--B model. Comm. Math. Sci. 15 (2017), 1265--1323}), where the existence of global-in-time finite-energy weak solutions was shown in two dimensional setting with stress diffusion. In the paper, we investigate the case without stress diffusion. We first restrict ourselves to the corotational setting as in \cite{LM00} ({\em P. L. Lions, N. Masmoudi. Global solutions for some Oldroyd models of non-Newtonian flows. Chin. Ann. Math., Ser. B 21(2) (2000), 131--146.}) We further assume the extra stress tensor is a scalar matrix and we derive a simplified model which takes a similar form as the multi-component compressible Navier--Stokes equations, where, however, the pressure term related to the scalar extra stress tensor has the opposite sign. By employing the techniques developed in \cite{VWY19, Novotny-Pokorny-19}, we can still prove the global-in-time existence of finite energy weak solutions in two or three dimensions, without the presence of stress diffusion.

\end{abstract}





\section{Introduction}

The Oldroyd--B model is a widely used constitutive model to describe the flow of viscoelastic fluids.  Different variants of such models were studied from different points of view within many decades: in the past mostly for incompressible fluids, recently, however, also for the compressible ones. There are many results in the context of small data problems, on the other hand, the number of existence results for large data without any restriction on the size of the data  or the length of the time interval are rather rare, even in the case of incompressible fluids. A typical feature for the viscoelastic fluid is the presence of an extra stress tensor which fulfils a certain type of transport equation. In most cases, the global-in-time existence results rely on the fact that an additional term describing the stress diffusion is present in these transport equations. Even though it is possible to justify the presence of the stress diffusion, in modelling, it is often neglected, as typically, such terms are many orders lower than other terms in the equations.  Our aim in the present paper  is to concentrate on a special case, where it is possible to neglect the stress diffusion in the compressible model and we can still obtain global-in-time existence of a solution without any restriction on the size of the data.

It is known that from the incompressible Navier--Stokes--Fokker--Planck system which is a micro-macro model describing incompressible dilute polymeric fluids one can derive the incompressible Oldroyd--B model in dumbbell Hookean setting, see \cite{Bris-Lelievre}. A similar derivation can be performed in the compressible setting, see \cite{Barrett-Lu-Suli}, 
where the existence of global-in-time finite-energy weak solutions was also shown in the two dimensional setting. However, an important role in the analysis of these models was played by the presence of the stress diffusion.

Let $\Omega \subset \mathbb{R}^d$ be a bounded open domain with a $C^{2,\beta}$ boundary (briefly, a $C^{2,\beta}$ domain), with $\beta \in (0,1]$, and $d = 2,3$. The compressible Oldroyd--B model derived in \cite{Barrett-Lu-Suli} posed in the time-space cylinder $Q_T:=(0,T)\times \O$ is the following:
\begin{alignat}{2}
\label{01a}
\d_t \vr + \Div_x (\vr \vu) &= 0,
\\
\label{02a}
\d_t (\vr\vu)+ \Div_x (\vr \vu \otimes \vu) +\nabla_x p(\vr)  - \Div_x \SSS(\nabla_x \vu) &=\Div_x \big(\TT - (kL\eta + \de\, \eta^2)\,\II\, \big)  +  \vr\, \ff,
\\
\label{03a}
\d_t \eta + \Div_x (\eta \vu) &= \e \Delta_x \eta,
\\
\label{04a}
\d_t \TT + {\rm Div}_x (\vu\,\TT) - \left(\nabla_x \vu \,\TT + \TT\, \nabla_x^{\rm T} \vu \right) &= \e \Delta_x \TT + \frac{k}{2\lambda}\eta  \,\II - \frac{1}{2\lambda} \TT.
\end{alignat}
Above, $\SSS(\nabla_x \vu)$ is the {\em Newtonian stress tensor} defined by
\be\label{Newtonian-tensor}
\SSS(\nabla_x \vu) = \mu^S \left( \frac{\nabla_x \vu + \nabla^{\rm T}_x \vu}{2} - \frac{1}{d} (\Div_x \vu) \II \right) + \mu^B (\Div_x \vu) \II,
\ee
where $\mu^S>0$ and $\mu^B\geq 0$ are the shear and bulk viscosity coefficients, respectively. The pressure $p$ and the density $\vr$ of the solvent are supposed to be related by the typical power law relation:
\be\label{pressure}
 p(\vr)=a \vr^\gamma, \quad a>0, \ \gamma \geq 1.
\ee
%

%
%

The extra stress tensor $\TT = (\TT_{ij})$, $1\leq i,j \leq d$ is a positive definite symmetric matrix defined on $Q_T$, and the notation ${\rm Div}_x(\vu\,\TT)$ is to be understood as
\be\label{def-Div-tau}
\left({\rm Div}_x(\vu\,\TT)\right)_{ij} = \Div_{ x}(\vu\,\TT_{ij}), \quad 1\leq i,j \leq d.\nn
\ee

The polymer number density $\eta$ is a non-negative scalar function defined as the integral of the probability density function $\psi$ in the conformation vector. The probability density function $\psi$ is governed by the Fokker--Planck equation. The conformation vector is a microscopic variable in the modelling of dilute polymer chains. The term $kL \eta + \de \eta^2$ in the momentum equation \eqref{02a} can be seen as the \emph{polymer pressure}, compared to the fluid pressure $p(\vr)$.

The meanings of the various quantities and parameters appearing in \eqref{01a}--\eqref{04a} were introduced in the derivation of the model in \cite{Barrett-Lu-Suli}.  In particular, the parameters $\e$, $k$, $\l$, $\de$, $L$ are all non-negative numbers.

\medskip

The mathematical study of Oldroyd--B models attracts a lot interests.  While, many fundamental problems are still open. Concerning the incompressible case without stress diffusion, the global-in-time existence of large data solutions is not known, even weak ones, either in the two dimensional or the three dimensional setting. With stress diffusion, the global-in-time existence of large data solutions in the two dimensional setting is known:  see \cite{Barrett-Boyaval} for weak solutions,  and \cite{CK} for strong solutions. But in the three dimensional setting, the global-in-time existence of large data solutions, strong or weak, is still open.  Note in this context that in \cite{KrPoSa} a special non-linear diffusion was used and global-in-time large-data solution for the corresponding Oldroyd--B model was possible to prove for a large variety of power-law models for the solvent stress tensor (from shear thinning to shear thickening, including, indeed, the linear dependence). Another interesting recent result is paper \cite{Bu_et_al}, where a slightly modified model was considered. The model was developed as a simplification of the general model based on the approach from \cite{RaSr} and global-in-time existence of weak solutions was shown for large data. 

\smallskip

 Even less is known concerning the compressible Oldroyd--B models. We recall some mathematical results for compressible viscoelastic models, which have been the subject of active research in recent years. First, note that in \cite{BuFeMa} a model based on the general approach from \cite{RaSr} was studied and global in time existence of weak solutions for sufficiently large $\gamma$ in the pressure law was shown. Note that the stress diffusion is present in the model. The existence and uniqueness of local strong solutions and the existence of global solutions near an equilibrium for macroscopic models of three-dimensional compressible viscoelastic fluids was considered in \cite{Qian-Zhang, Hu-Wang3, Hu-Wu, Lei}. In particular, Fang and Zi \cite{Fang-Zi} proved the existence of a unique local-in-time strong solution to a compressible Oldroyd--B model and established a blow-up criterion for strong solutions.  In \cite{Barrett-Lu-Suli}, not only the compressible Oldroyd--B model \eqref{01a}--\eqref{pressure} was derived, but also the existence of global-in-time weak solutions in two dimensional setting was shown. Recently in \cite{Lu-Z.Zhang18}, one of the authors and his collaborator proved the weak-strong uniqueness, gave a refined blow-up criterion and showed a conditional regularity result in two dimensional setting. 

\medskip

There are stress diffusion terms $\e \Delta_x \eta$ and $\e \Delta_x \TT$ in \eqref{03a} and \eqref{04a} which help for the mathematical analysis of existence theory. Such spatial stress diffusions are allowed in some models of complex fluids, such as the creeping flow regime, as pointed out in \cite{CK}. Also in the modelling of the compressible Navier--Stokes--Fokker--Planck system arising in the kinetic theory of dilute polymeric fluids, where polymer chains immersed in a barotropic, compressible, isothermal, viscous Newtonian solvent, Barrett and S\"uli \cite{Barrett-Suli, Barrett-Suli1, Barrett-Suli2, Barrett-Suli4, BS2016} observed the presence of the centre-of-mass diffusion term $\e \Delta_x \psi$, where $\psi$ is the probability density function depending on both microscopic and macroscopic variables; as a result, its macroscopic closure (the compressible Oldroyd--B model) contains such diffusion terms.  The center-of-mass coefficient $\e = (l_0^2/ L_0^2) ({1}/{4(K+1)\l})$, where $L_0$ is the macroscopic characteristic length-scale of the solvent flow and $l_0$ is the characteristic microscopic length-scale. The parameter $\l>0$ is the Deborah number, and $K+1$ is the number of beads in the bead-spring polymer chain.

However, in standard derivations of bead-spring models the center-of-mass diffusion term is routinely omitted, on the grounds that it is several orders of magnitude smaller than the other terms in the equations. Indeed, Bhave, Armstrong and Brown \cite{BAB91} show the ratio $l^2_0/L_0^2$ to be in the range of about $10^{-9}$ to $10^{-7}$. For such a reason, in most previous mathematical studies of Oldroyd--B model, the stress diffusion is not included, for example in  Renardy \cite{Renardy90}, Guillop\'e and Saut \cite{G-S90-1, G-S90-2} and  Fern\'andez-Cara, Guill\'en and Ortega \cite{F-G-O02} for the local-in-time well-posedness, as well as the global-in-time well-posedness with small data, Lions and Masmoudi \cite{LM00} for the global-in-time existence of weak solutions with large initial data in the corotational derivative setting. These mentioned results considered incompressible Oldroyd-B models. While for the compressible models without stress diffusion, according to the authors' knowledge, there is basically no results concerning the global-in-time existence of large solutions. This is the direction that we are working on in this paper and the aim is to find out under which further assumptions the global-in-time existence of large solutions can be shown for the compressible Oldroyd-B model without stress diffusion.

Inspired by the work of Lions and Masmoudi \cite{LM00}, we consider the corotational derivative setting. Additionally, assuming the extra stress tensor is a scalar matrix, we derive a simplified model which takes a similar form as the multi-component compressible Navier--Stokes equations. However, as explained in the next section, the simplified model we derive has a pressure term that has a {\em wrong} sign which could cause the a priori energy estimates fail. While under a domination assumption on the data, we can still employ the techniques recently developed in \cite{VWY19, Novotny-Pokorny-19} and prove the global-in-time existence of large date finite energy weak solutions in two or three dimensions.

\section{Formulation of the problem and main results}

In this section, we will formulate our simplified models step by step and give the main results. As mentioned in the introduction, we will not include the stress diffusion $\e \Delta_x \eta$ and $\e \Delta_x \TT$ in \eqref{03a}--\eqref{04a}.

\subsection{Formulation of the problem--a simplified model}

In \eqref{04a}, the derivative
\be\label{upper-derivative}
\d_t \TT + {\rm Div}_x (\vu\,\TT) - \left(\nabla_x \vu \,\TT + \TT\, \nabla_x^{\rm T} \vu \right)
\ee
is called the upper convected derivative which naturally appears in several macroscopic models derived from microscopic models. It is known to be frame invariant. Other frame
invariant derivatives exist, like the corotational one (see \cite{Bris-Lelievre,LM00})
\be\label{upper-derivative-corotational}
\d_t \TT + {\rm Div}_x (\vu\,\TT) - \left(\o( \vu) \,\TT - \TT\,  \o( \vu) \right),
\ee
where $\o(\vu) = \big(\nabla_x \vu - \nabla_x^{\rm T} \vu \big)/2$ is the vorticity tensor. We first restrict ourselves to the corotational derivative setting \eqref{upper-derivative-corotational} in \eqref{04a}, where the deformation tensor ${\mathbb{D}}(\vu) = \big(\nabla_x \vu - \nabla_x^{\rm T} \vu \big)/2$ is neglected in the upper convected derivative \eqref{upper-derivative}.

\medskip

We further make a serious simplifying assumption that the extra stress tensor $\TT$ is a scalar matrix:
\be\label{TT-scalar-ass}
\TT = \tau \II, \quad \mbox{for some scalar function $\tau$.}
\ee
Since $\TT$ is a positive definite matrix, $\tau$ is supposed to be a positive scalar function. Then equation \eqref{04a} without stress diffusion $\e \Delta_x \TT$ becomes
\be\label{04aa}
\d_t \tau + \dive_x (\tau\vu)  =  \frac{k}{2\lambda}\eta  - \frac{1}{2\lambda} \tau.
\ee
By introducing
\be\label{tilde-tau}
\tilde \tau = \tau  - k \eta,
\ee
we deduce from \eqref{03a} and \eqref{04aa} that
\be\label{04aaa}
\d_t \tilde \tau + \dive_x (\tilde \tau \vu)  = - \frac{1}{2\lambda} \tilde \tau.
\ee

\medskip

By omitting the tilde in \eqref{04aaa}, and collecting \eqref{01a}, \eqref{02a} and \eqref{03a} without stress diffusion, we finally derive the following model in $Q_T = (0,T)\times \O$:
\begin{alignat}{2}
\label{01}
\d_t \vr + \Div_x (\vr \vu) &= 0,
\\
\label{02}
\d_t (\vr\vu)+ \Div_x (\vr \vu \otimes \vu) +\nabla_x \big(p(\vr) + q(\eta) -\tau\big) - \Div_x \SSS(\nabla_x \vu) &=  \vr\, \ff,
\\
\label{03}
\d_t \eta + \Div_x (\eta \vu) &= 0,
\\
\label{04}
\d_t  \tau + \dive_x ( \tau \vu)  &= - \frac{1}{2\lambda} \tau.
\end{alignat}
Here the pressure $p(\vr)$ and the Newtonian stress tensor $\SSS(\nabla_x \vu)$ are defined as in \eqref{pressure} and \eqref{Newtonian-tensor}, and the polymer pressure $q(\eta)$ takes the form
\be\label{q-eta}
q(\eta): = k(L-1) \eta + \de\,\eta^2.
\ee
We impose for simplicity the no slip boundary condition
\be\label{boundary}
\vu = \vc{0} \quad \mbox{on $(0,T)\times \d \O $}.
\ee
The external force $\ff$ is assumed to be in $L^\infty(Q_T;\R^d)$. In this paper, we focus on the simplified model \eqref{01}--\eqref{boundary} and we will show the global-in-time existence of large data finite energy weak solutions. Note that it is natural to assume that $\varrho$ and $\eta$ are non-negative, while the fact that the original $\tau$ is non-negative does not say anything about the sign of $\tilde \tau$ defined in \eqref{tilde-tau}. However, in what follows, we will assume that $\tilde \tau$, i.e. our $\tau$ in system \eqref{01}--\eqref{04}, is non-negative.

\medskip

Having a first glance at this model, it looks like the multi-component compressible Navier--Stokes equations considered in \cite{VWY19,Novotny-Pokorny-19}, with three density functions $\vr, \eta, \tau$. But the sign of the pressure term related to $\tau$ in the momentum equations is negative, which is opposite to the sign of the other two pressure terms $p(\vr)$ and $q(\eta)$. This different sign is {\em bad} or {\em wrong} in the sense that in the a priori energy equality, the energy functional could have no determinate sign, even after some modifications by adding a fixed large constant. This can be seen later in the energy inequality \eqref{energy1-f}.

Without a clear sign of the energy functional in the energy inequality one could fail to deduce uniform a priori estimates on the solutions which is usually the starting point of the existence theory. We overcome this difficulty by imposing a domination assumption on the initial data. One can show that this domination will preserve as time goes. This allows us to obtain uniform estimates on the solutions from the energy inequality with {\em bad} sign.

Another difference compared to the multi-component compressible Navier--Stokes equations considered in \cite{VWY19,Novotny-Pokorny-19} is that in the momentum equation, the time derivative term and the convective term involves only the fluid density $\vr.$ Fortunately, as we will see below, this does not cause any troubles. The a priori estimates are not really influenced. Moreover, by uniform estimates and Arzel\`a--Ascoli type lemma, we can show the weak compactness of the convective terms (i.e., terms $\vr \vu$ and $\vr\vu\otimes \vu$). Thus in the proof of the weak compactness of the pressures by using effective viscous flux, these two terms do not play any role. Then the strong convergence of the density can still be shown.

\subsection{Global-in-time finite energy weak solutions}

The basic hypotheses on the initial data are
\ba\label{ini-data-f}
&\vr(0,\cdot) = \vr_0(\cdot) \ \mbox{with}\ \vr_0 \geq 0 \ {\rm a.e.} \ \mbox{in} \ \O, \quad \vr_0 \in L^\gamma(\O),
\\
&\vr_0\vu_0  \in L^1(\O;\R^d),\ \  \vr_0|\vu_0|^2 \in L^1(\O), \quad \vu_0 = \vu(0,\cdot),\\
&\eta(0,\cdot)=\eta_0(\cdot) \ \mbox{with}\ \eta_0 \geq 0 \ {\rm a.e.} \ \mbox{in} \ \O, \quad \eta_0 \in L^2 (\O), \\
&\tau(0,\cdot) = \tau_0(\cdot) \ \mbox{with}\ \tau_0 \geq 0 \ \mbox{a.e. in}  \ \O,  \quad  \tau_0 \log \tau_0 \in L^1(\O).
\ea
We give the definition of finite energy weak solutions:
\begin{definition}\label{def-weaksl-f} Let $T>0$ and $\O\subset \R^d$ be a bounded $C^{2,\beta}$ domain with $0<\beta\leq 1$. We say that $(\vr,\vu,\eta,\tau)$ is a finite-energy
weak solution in $Q_T$ to the system of equations \eqref{01}--\eqref{boundary},  supplemented by the initial data \eqref{ini-data-f}, if:
\begin{itemize}
\item $\vr \geq 0 \ {\rm a.e.\  in} \ (0,T) \times \Omega$, $\vr \in  C_w([0,T];  L^\gamma(\Omega))$,  $\vu\in L^{2}(0,T;W_0^{1,2}(\Omega; \R^d))$,
\vspace{-2mm}
\begin{align*}
&\vr \vu \in C_w([0,T]; C^{1}(\Omega; \R^d)),\quad \vr |\vu|^2 \in L^\infty(0,T; L^{1}(\Omega)),\\
&\eta  \geq 0  \ {\rm a.e.\  in} \ (0,T) \times \Omega,\quad \eta \in C_w ([0,T]; L^2(\Omega)), \\
&\tau \geq 0  \ {\rm a.e.\  in} \ (0,T) \times \Omega, \quad \tau\log\tau \in C_w([0,T]; L^1(\Omega)).
\end{align*}
\item For any $t \in (0,T)$ and any test function $\phi \in C^\infty([0,T] \times \Ov{\Omega})$, one has
\be\label{weak-form1-f}
\int_0^t\intO{\big[ \vr \partial_t \phi + \vr \vu \cdot \Grad \phi \big]} \,\dt' =
\intO{\vr(t, \cdot) \phi (t, \cdot) } - \intO{ \vr_{0} \phi (0, \cdot) },
\ee
\be\label{weak-form2-f}
\int_0^t \intO{ \big[ \eta \partial_t \phi + \eta \vu \cdot \Grad \phi\big] } \, \dt' =  \intO{ \eta(t, \cdot) \phi (t, \cdot) } - \intO{ \eta_{0} \phi (0, \cdot) },
\ee
\be\label{weak-form3-f}
\int_0^t \intO{ \big[ \tau \partial_t \phi + \tau \vu \cdot \Grad \phi - \frac{1}{2\l} \tau \phi\big]} \, \dt' =  \intO{ \eta(t, \cdot) \phi (t, \cdot) } - \intO{ \eta_{0} \phi (0, \cdot) }.
\ee

\item For any $t \in (0,T)$ and any test function $\vvarphi \in C^\infty([0,T]; \DC({\Omega};\R^d))$, one has
\ba\label{weak-form4-f}
&\hspace{-0.7cm}\int_0^t \int_\O \big[ \vr \vu \cdot \partial_t \vvarphi + (\vr \vu \otimes \vu) : \Grad \vvarphi  + \big(p(\vr)+q(\eta)  - \tau\big) \Div_x \vvarphi  - \SSS(\nabla_x \vu) : \Grad \vvarphi \big] \dx \, \dt'\\
&\hspace{-0.7cm}= - \int_0^t \intO{   \vr\, \ff \cdot \vvarphi } \, \dt' + \intO{ \vr \vu (t, \cdot) \cdot \vvarphi (t, \cdot) } - \intO{ \vr_{0} \vu_{0} \cdot \vvarphi(0, \cdot) }.
\ea
%

%
%
%
%
%
%
%

%
%
\item  For a.e. $t \in (0,T)$, the following \emph{energy inequality} holds
\ba\label{energy1-f}
&\int_\O \left[ \frac{1}{2} \vr |\vu|^2  + H(\vr,\eta, \tau) \right](t,\cdot)\,\dx  +  \int_0^t  \int_\O \SSS(\nabla_x \vu): \nabla_x \vu\,\dx \,\dt'  \\
& \leq  \int_\O \left[ \frac{1}{2} \vr_0 |\vu_0|^2 + H(\vr_0,\eta_0, \tau_0)  \right]\dx + \int_0^t \int_\O \vr\,\ff \cdot  \vu \,\dx\,\dt' + \frac{1}{2\lambda} \int_0^t \int_\O (\tau \log \tau + \tau)  \,\dx \,\dt',
\ea
where the Helmholtz free energy is defined as
\be\label{pressure-entropy}
H(\vr, \eta, \tau) = P(\vr) + Q(\eta) -   \tau \log \tau
\ee
with
\be\label{pressure-entropy-2}
Q(\eta) =  \de \,\eta^2 + k (L-1)  \eta \log \eta , \quad P(\vr) = \left\{ \begin{aligned} & \frac{a}{\g-1} \vr^\g, \ \ &&\mbox{if $\g\neq 1$},\\
&  a \vr\log \vr, \ \ && \mbox{if $\g = 1$}.
\end{aligned} \right.
\ee
\end{itemize}
\end{definition}

We are now in position to state our main result.
\begin{theorem}\label{thm}
Let $\Omega \subset \mathbb{R}^d, \, d=2,3$, be a bounded $C^{2,\beta}$ domain with $\beta \in (0,1]$. Let $0 < \g\leq 2$, the constant parameters $\l$, $\de$ be positive, and $k, L$ be non-negative. We further assume that the initial data satisfy the domination relation:
\be\label{ini-domination}
 \vr_0 \leq \overline C \eta_0, \quad  \tau_0 \leq \overline C \eta_0 \quad \mbox{a.e. in } \Omega \mbox{ for some $\overline C > 0$}.
\ee
Then, for any $T>0$, there exists a finite-energy weak solution $(\vr,\vu,\eta,\tau)$ in the sense of Definition \ref{def-weaksl-f} with initial data \eqref{ini-data-f} by replacing the integrability on $\vr$ and $\tau$ by
\be
\vr\in C_w ([0,T];  L^2(\Omega)), \quad \tau \in C_w ([0,T];  L^2(\Omega)).\nn
\ee
Moreover, the domination condition preserves for all times:
\be\label{ini-domination-2}
 \vr(t,x) \leq \overline C \eta(t,x), \quad   \tau(t,x) \leq \overline C \eta(t,x) \quad \mbox{for a.a. $(t,x)\in Q_T$}.
\ee
\end{theorem}

\begin{remarks}\label{remarks-thm} Before presenting the proof, we give several remarks concerning our result:

\begin{itemize}

\item We first remark that $\g\leq 2$ is not essential. By \eqref{ini-data-f} and \eqref{ini-domination}, one has $\vr_0\in L^2(\O)$, $\tau_0 \in L^2(\O)$. If $\g>2$, $\vr_0 \in L^\g(\O)$ has better integrability than $L^2$. This case is actually easier: one can use $\vr$ as the benchmark density and consider the following domination condition:
\be\label{ini-domination-new}
\eta_0 \leq \overline C \vr_0, \quad  \tau_0 \leq \overline C \vr_0 \quad \mbox{for some $\overline C >  0$}.\nn
\ee
Then the results and the proofs follow in the same manner. Another possibility is to employ the idea in \cite{Wen19} to avoid any domination conditions. For such a reason, we will consider only the case $\g \leq 2$.

\item In general, the adiabatic number $\g \geq1$. As we see in the energy inequality \eqref{energy1-f} and the domination \eqref{ini-domination-2} (also in other aspects, for example in the derivation of higher integrability of the pressure), the larger $\gamma$ is, the better integrability one has for density $\vr$, and the better result one can get. Due to the domination relations in \eqref{ini-domination} and the quadratic term $\eta^2$ in $q(\eta)$, by using the argument in \cite{Novotny-Pokorny-19}, we can relax the restriction on $\g$ to any $\g>0$.

\item Assumption \eqref{ini-domination} on the initial data means that one density dominates the others. This allows us to control $\vr$ and $\tau$ by using the estimates on $\eta$. This assumption is inspired by the study of global-in-time existence of weak solutions for multi-component fluid flows, see \cite{VWY19} and \cite{Novotny-Pokorny-19}. Recently in \cite{Wen19} the author eliminates this condition under some further restrictions on the adiabatic numbers. Note that it would lead to restriction $\g \geq \frac 95$.

\item Other boundary conditions such as periodic boundary condition and Navier boundary condition can also be considered. The results and the proofs follow in a straightforward way as the present case with homogeneous Dirichlet boundary condition under proper modifications.

\item The $C^{2,\b}$ regularity assumption on the domain $\O$ is not necessary and could be relaxed. For example, Lipschitz regularity for bounded domain $\O$ will be enough to prove Theorem \ref{thm} via a smoothing technique on the boundary, see \cite{FNP-02}.

\end{itemize}
\end{remarks}

\subsection{Preliminaries}

Before giving the proof, we recall the results from \cite{Novotny-Pokorny-19} with a few comments on the hypotheses. Therein, the following problem is studied

\begin{equation}\label{eq8.1}
\begin{aligned}
\partial_t \vr  &+ \Div(\vr \vu) = 0, \\
\partial_t Z_i &+ \Div (Z_i \vu)  = 0,  \quad i=1,2,\dots, K,\\
\partial_t \big((\vr+\sum_{i=1}^K Z_i)\vu\big) &+ \Div\big((\vr+\sum_{i=1}^K Z_i) \vu\otimes \vu) + \nabla P(\vr,Z_0,Z_1,\dots, Z_K)
=  \mu \Delta \vu + (\mu+\lambda)\Grad \Div \vu,
\end{aligned}
\end{equation}
together with the boundary condition
$
\vu = \vc{0}
$
on $(0,T)\times \partial \Omega$, and the initial conditions in $\Omega$
\begin{equation} \label{eq8.3}
\begin{aligned}
\vr(0,x) &= \vr_0(x), \\
Z_i(0,x) &= Z_{i0}(x), \quad i=1,2,\dots, K,\\
\Big(\vr + \sum_{i=1}^K Z_i\Big)\vu(0,x) &= \vc{m}_0(x).
\end{aligned}
\end{equation}

The weak formulation of this problem is similar to our weak formulation from Definition \ref{def-weaksl-f}.

As we will explain later, the fact that for our problem one of the continuity equations has non-trivial special right hand side does not play any significant role in the analysis. Similarly, the fact that in our problem we have only $\partial_t(\vr \vu)$ and $\Div(\vr \vu \otimes \vu)$ rather simplifies the proof. Also the presence of the external force $\vr \vc{f}$ with $\vc {f}\in L^\infty(Q_T;\R^3)$ does not cause any troubles. Finally, the 2-D case is even simpler and the hypotheses presented below are surely sufficient to get a solution also in this situation.
\\

\noindent {\bf Hypothesis (H1).}
\begin{equation} \label{eq2.1MF}
\begin{aligned}
(\vr_0,Z_{10}, Z_{20}, \dots, Z_{K0})\in  {\cal O}_{\underline {\vec {a}}} := \big\{(\vr,Z_1,Z_2,\dots, Z_K)\in \R^{K+1} |\vr\in [0,\infty), \,\underline {a}_i \vr < Z_i < \overline {a}_i \vr \big\},
\end{aligned}
\end{equation}
where $0\leq \underline{a}_i<\overline{a}_i <\infty$, $i=1,2,\dots, K$.
\\

\noindent {\bf Hypothesis (H2).}
\ba \label{eq2.6MF}
&\vr_0 \in L^\gamma(\Omega), \; Z_{i0} \in L^{\beta_i}(\Omega) \; \text{ if } \beta_i > \gamma,\\
&\vc{m}_0 \in L^1(\Omega;\R^3),\; (\vr_0+\sum_{i=1}^K Z_{i0})|\vu_0|^2\in L^1(\Omega), \, i=1,2,\dots, K.
\ea
\\

\noindent {\bf Hypothesis (H3).}

Function $P\in C(\overline{{\cal O}_{\underline {\vec{a}}}})\cap C^1({\cal O}_{\underline{\vec{ a}}})$
and
\begin{equation}\label{?!+++}
\forall\vr\in (0,1),\; \sup_{s\in { \Pi_{i=1}^K}[{\underline a}_i,{\overline a}_i]} |P(\vr,\vr s_1,\vr s_2,\dots, \vr s_K)|\le C\vr^\alpha\;\mbox{with some $ C>0$ and $\alpha>0$},
\end{equation}
and
\begin{equation} \label{eq8.7}
\underline{C}(\vr^\gamma + \sum_{i=1}^K Z_i^{\beta_i} -1) \leq P(\vr,Z_1,\dots, Z_K) \leq \overline{C}(\vr^\gamma + \sum_{i=1}^K Z_i^{\beta_i} +1)
\;\mbox{in ${\cal O}_{\underline {\vec{a}}}$}
\end{equation}
with $\gamma \geq \frac 95$, $\beta_i >0$, $i=1,2,\dots, K$.
We moreover assume for $i=1,2,\dots, K$
\begin{equation}\label{MFeq2.5-}
|\partial_{Z_i}P(\vr,Z_1,Z_2,\dots, Z_K)|\le C(\vr^{-\underline\Gamma}+\vr^{\overline\Gamma-1})\;\mbox{in ${\cal O}_{\underline {\vec{a}}}$}
\end{equation}
with some $0\le\underline\Gamma<1$, and with some $0< \overline\Gamma < \gamma + \gamma_{BOG}$ if $\underline{a}_i=0$,
$0<\overline\Gamma<{\rm max}\{\gamma+\gamma_{BOG}, \beta_i+(\beta_i)_{BOG}\} $ if $\underline{a}_i>0$. \\

\noindent {\bf Hypothesis (H4).}
We assume
%
\begin{equation} \label{MFeq2.4}
P(\vr,\vr s_1,\vr s_2\dots, \vr s_K)=
{\cal P}(\vr,s_1,s_2,\dots, s_K) - {\cal R} (\vr,s_1,s_2,\dots, s_K),
\end{equation}
where $[0,\infty)\ni\vr\mapsto {\cal P}(\vr,s_1,s_2,\dots, s_K)$ is non decreasing  for any $s_i\in [\underline {a}_i,\overline {a}_i]$, $i=1,2,\dots, K$, and $\vr\mapsto {\cal R}(\vr,s_1,s_2,\dots, s_K)$ is for any $s_i \in [\underline {a}_i,\overline {a}_i]$, $i=1,2,\dots, K$ a non-negative $C^2$-function  in $[0,\infty)$ with uniformly bounded $C^2$-norm with respect to $s_i \in [\underline {a}_i, \overline {a}_i]$, $i=1,2,\dots, K$ and with compact support uniform with respect to $s_i \in [\underline {a}_i, \overline {a}_i]$, $i=1,2,\dots, K$. Here, $\underline {a}_i, \overline {a}_i$ are the constants from relation (\ref{eq2.1MF}). \\

\noindent {\bf Hypothesis (H5).}
Functions $\vr\mapsto P(\vr,Z_1,Z_2,\dots,Z_K)$, $Z_i>0$, $i=1,2,\dots, K$ resp. \linebreak $(Z_1,Z_2,\dots, Z_K)\mapsto \partial_{Z_j} P(\vr,Z_1,Z_2,\dots, Z_K)$, $\vr>0$,
are Lipschitz on $\cap_{i=1}^K(Z_i/{{\overline a}_i},Z_i/{{\underline a}_i})\cap(\underline r,\infty)^K$ resp.
$\Pi_{i=1}^K({{\underline a}_i}\vr, \overline{a}_i\vr)\cap(\underline r,\infty)^K$ for all $\underline r>0$  with Lipschitz constants
\begin{equation}\label{eq2.3a-MF}
\widetilde L_P\le C(\underline r)(1+|Z|^A)\;\mbox{ resp.}\; \widetilde L_P\le C(\underline r)(1+\vr^A)
\end{equation}
with some non negative number $A$. Number $C(\underline r)$ may diverge to $+\infty$ as $\underline r\to 0^+$.
%
\\

The following result is taken from \cite[Theorem 15]{Novotny-Pokorny-19}.

\begin{theorem} \label{t2}
Let $\gamma > \frac 95$. Then under Hypotheses (H1--H5), there exists at least one weak solution to problem \eqref{eq8.1}--\eqref{eq8.3}. Moreover, the densities $\vr \in C_{weak}([0,T); L^\gamma (\Omega))$, \linebreak $Z_i \in C_{weak}([0,T); L^{\max\{\gamma,\beta_i\}} (\Omega))$, $i=1,2,\dots, K$, $(\vr+\sum_{i=1}^K Z_i)\vu \in C_{weak}([0,T); L^q (\Omega;\R^3))$ for some $q>1$, and $P(\vr,Z_1,Z_2,\dots, Z_K) \in L^q(\Omega)$ for some $q>1$.
\end{theorem}

Note that in \cite{Novotny-Pokorny-19}, also the case $\gamma=\frac 95$ is treated. Since it requires certain extra conditions which we do not need (recall that our largest exponent is equal to $2$), we skip them.

The above Hypotheses, in particular (H3) and (H5), are connected with the Helmholtz free energy $H_P(\vr,Z_1,\dots, Z_K)$, a solution to the partial differential equation of the first order in ${\cal O}_{\underline{\vec{ a}}}$,
\ba\label{MFeqH}
P(\vr,Z_1,Z_2,\dots, Z_K) &= \vr \frac{\partial H_P(\vr,Z_1,Z_2,\dots, Z_K)}{\partial \vr} \\
&\quad +  \sum_{i=1}^K Z_i\frac{\partial H_P(\vr,Z_1,Z_2,\dots, Z_K)}{\partial Z_i}-H_P(\vr,Z_1,Z_2,\dots, Z_K)
\ea
in the form
\begin{equation} \label{MFsolH}
H_P(\vr,Z_1,Z_2,\dots, Z_K) = \vr \int_{1}^{\vr} \frac{P\big(s, s\frac{Z_1}{\vr}, s\frac{Z_2}{\vr}, \dots, s\frac{Z_K}{\vr}\big)}{s^2} \, {\rm d}s,\;
H_P(0,\ldots,0)=0.
\end{equation}

However, we consider a slightly different form of the solution to \eqref{MFeqH}, using the fact that the pressure can be written as a sum of three pressures, each dependent only on one unknown. Therefore we have to modify the Hypotheses (H3) and (H5) for this situation. We will comment on this in the next section.

We also formulate one important auxiliary result which is the main ingredient of the compactness of the densities other than {\it the main one}, see \cite[Proposition 7]{Novotny-Pokorny-19}.

\begin{proposition}\label{P2.5}
\begin{description}
\item{\it 1.}
Let
$$
\vu_n\in L^2(I,W_0^{1,2}(\Omega;\R^3)),\;
(\vr_n, Z_n)\in {\cal O}_0\cap \Big(C(\overline I;L^1(\Omega))\cap L^2(Q_T)\Big)^2.
$$
Suppose that
$$
\sup_{n\in N}\big(\|\vr_n\|_{L^\infty(I;L^\gamma(\Omega))}+\|Z_n\|_{L^\infty
(I;L^{\gamma}(\Omega))} + \|\vr_n\|_{L^2(Q_T)}+ \|\vu_n\|_{L^2(I;W^{1,2}(\Omega))}\big) <\infty,
$$
where $\gamma>\frac 65$,
and that both couples $(\vr_n,\vu_n)$, $(Z_n,\vu_n)$ satisfy continuity equation
$$
\partial_t \vr_n + \Div(\vr_n\vu_n) = 0, \qquad \partial_t Z_n + \Div(Z_n\vu_n) = 0.
$$
Then,
up to a subsequence (not relabeled)
$$
\begin{aligned}
(\vr_n, Z_n)\to (\vr,Z)\;\mbox{in $(C_{\rm weak}(\overline I;L^\gamma(\Omega)))^2$}, \quad\vu_n\rightharpoonup\vu\;\mbox{weakly in $L^2(I;W^{1,2}(\Omega;\R^3))$,}
\end{aligned}
$$
where $(\vr,Z)$ belongs to spaces
$$
{\cal O}_0\cap (L^2(Q_T))^2\cap (L^\infty(I,L^\gamma(I,\Omega)))^2\cap (C(\overline I;L^1(\Omega))^2
$$
and $(\vr,\vu)$ as well as $(Z,\vu)$ verify continuity equation in the renormalized sense.
\item{\it 2.} We define  for all $t\in \overline I$,
\begin{equation}\label{sn}
s_n(t,x)=\frac{Z_n(t,x)}{\vr_n(t,x)},\quad s(t,x)=\frac {Z(t,x)}{\vr(t,x)},\nn
\end{equation}
where $s_n(t,x) =0$ ($s(t,x) =0$) provided $\vr_n(t,x) =0$ ($\vr(t,x) =0$).
Suppose in addition to assumptions of Item 1. that
$$
\int_{\Omega}{\vr_n(0,x)}s_n^2(0,x){\rm d} x\to \int_{\Omega}\vr(0,x) s^2(0,x)\,{\rm d} x.
$$
Then $s_n, s\in C(\overline I;L^q(\Omega))$, $1\le q<\infty$ and  for all $t\in \overline I$,
$0\le s_n(t,x)\le\overline a$, $0\le s(t,x)\le \overline a$ for a.a. $x\in \Omega$. Moreover, both
$(s_n,\vu_n)$ and $(s,\vu)$ satisfy
transport equation
$$
\partial_t s_n + \vu_n \cdot \nabla s_n = 0, \qquad \partial_t s + \vu \cdot \nabla s = 0
$$
in the weak and time-integrated sense (cf. \eqref{renormal-3} below).
\item{\it 3.}
Finally,
\begin{equation}\label{cvs}
\int_{\Omega}(\vr_n|s_n-s|^\th)(\tau,\cdot)\,{\rm d} x\to 0\;\mbox{with any $1\le \th<\infty$}\nn
\end{equation}
for all $\tau\in [0,T]$.
\end{description}
\end{proposition}

The rest of the paper is devoted to the proof of Theorem \ref{thm}. First we check that we fulfil most of the Hypotheses (H1--H5) presented above and then explain that the remaining ones are actually not important. To document this fact, we also give the most important ideas of the proof: the a priori estimates and the weak compactness of the solutions to our problem. The way of constructing approximate solution sequences by multiple levels of approximations is now well understood for compressible Navier--Stokes equations, see \cite{F-book,FNP,N-book, Novotny-Pokorny-19}, hence we will not repeat the approximation schemes. The full proof follows by combining the rather classical construction of approximate solutions done in the above references and the uniform estimates and the compactness shown in the following sections.

In the sequel, we use $C$ to denote a universal positive constant whose value may differ from line to line.

\section{Proof of the main result}

\subsection{Properties of the pressure}

The main purpose of this subsection is to investigate the properties of the pressure and verify Hypotheses (H1--H5) presented above. Before doing so, note that the role of the function $\vr$ in Theorem \ref{t2} is in our case of Theorem \ref{thm} played by the function $\eta$. Therefore the exponent $\gamma$ in Theorem \ref{t2} is equal to $2$ and $\beta_1 = \gamma$ ($Z_1 =\vr$) and $\beta_2 = 1$ ($Z_2=\tau$).

With this notation, it is not difficult to see that Hypothesis (H1) is fulfilled with $\underline{a}_1 =0$ and $\underline{a}_2 =0$.  Moreover, $\overline{a}_1=\overline{a}_2 =\overline{C}>0$. Next, Hypothesis (H2) is also fulfilled (however, with a slight straightforward  modification due to the different form of the momentum equation).

\medskip

We denote the total pressure as
\be\label{presure-total}
h(\eta,\vr,\tau) : =  q(\eta) + p(\vr) -\tau =  \de\,\eta^2 + k(L-1) \eta + a \vr^\g  - \tau.
\ee

Related to the domination condition \eqref{ini-domination}, we denote the set
\be\label{def-barS-1}
S: = \{(\eta,\vr,\tau) \in \R^3: 0 < \vr < \overline C \eta, \quad  0 < \tau < \overline C \eta\}
\ee
with closure
\be\label{def-barS-2}
\overline {S}: = \{(\eta,\vr,\tau) \in \R^3: 0  \leq \vr \leq \overline C \eta, \quad 0 \leq \tau \leq \overline C \eta\}.
\ee%
Then $S$ plays the role of ${\mathcal O}_{\vec{a}}$ and $\overline {S}$ of its closure in the above hypotheses.

Clearly, the total pressure $h(\eta,\vr,\tau) \in C(\overline {S})$ and $h(\eta,\vr,\tau) \in C^1(S)$. 
For all $\eta \in (0,1)$ and for all $(\eta,\vr,\tau)\in \overline {S}$, direct calculation gives
\ba\label{H3-1}
|h(\eta,\vr,\tau)| \leq \de \eta^2 + k |L-1| \eta + a \overline C^\g \eta^\g   + \overline C \eta \leq C (\eta+\eta^{\g}) \leq C\eta.
\nn
\ea
Next, in $S$, we have that
\be\label{press}
C_1 (\eta^2 +\vr^\gamma -\tau-1) \leq h(\eta,\vr,\tau) \leq C_2 (\eta^2 +\vr^\gamma +\tau+1)
\nn
\ee
for some positive constants $C_{1}, C_{2}$.
Using the domination assumption and the resulted domination for all times (proved, however, in the following subsection), we get that
$$
-\tau \geq -\overline{C}\eta;
$$
whence we have
\be\label{press_1}
C_1 (\eta^2 +\vr^\gamma -1) \leq h(\eta,\vr,\tau) \leq C_2 (\eta^2 +\vr^\gamma +\tau+1).\nn
\ee
Using once more the domination, we can also write
\be\label{press_2}
C_1 (\eta^2 +\vr^\gamma +\tau -1) \leq h(\eta,\vr,\tau) \leq C_2 (\eta^2 +\vr^\gamma +\tau+1).\nn
\ee
Moreover, it is not difficult to see that we also have
\be\label{Helm}
C_1 (\eta^2 +\vr^\gamma +\tau |\log \tau| -1) \leq H(\eta,\vr,\tau) \leq C_2 (\eta^2 +\vr^\gamma +\tau |\log \tau| +1),\nn
\ee
which follows from the form of our Helmholtz free energy and the domination properties. This estimate is important in the construction of weak solutions to our problem.
Furthermore, for each $(\eta,\vr,\tau)\in S$,
\ba\label{H3-2}
|\d_\tau h(\eta,\vr,\tau)| = 1, \quad |\d_\vr h(\eta,\vr,\tau)| = |a \g \vr^{\g-1}|.\nn
\ea
For $\g\geq 1$ it implies that \eqref{MFeq2.5-} is fulfilled for the choice $\underline{a}_1 =0$. However, for $\g \in (0,1)$ we cannot fulfil this assumption for $\underline{a}_1 =0$, as $\g-1<0$ and we need to control the function $\vr$ by $\eta$ from below. However, this condition is in fact in our case not needed and we have an alternative way how to overcome its use. It is connected with the proof that the pressure $h(\eta_n,\vr_n,\tau_n)$ converges weakly in $L^1(Q_T)$ to $h(\eta,\vr,\tau)$ and it will be explained in Subsection \ref{sec:str-cov-density}. Thus, the main part of Hypothesis (H3) is satisfied.

\medskip

For each $(\eta,\vr,\tau)\in \overline {S}$, we define the following functions
\be\label{def-s}
s_\vr : =\left\{ \begin{aligned} & \frac{\vr}{\eta}, \quad && \mbox{if $\eta>0$}, \\
& 0, \quad && \mbox{if $\eta=0$},
\end{aligned} \right.
\quad
s_\tau : =\left\{ \begin{aligned} & \frac{\tau}{\eta}, \quad && \mbox{if $\eta>0$}, \\
& 0, \quad && \mbox{if $\eta=0$}.
\end{aligned} \right.
\ee
Clearly $s_\vr, s_\tau \in [0, \overline C]$ for all $(\eta,\vr,\tau)\in \overline {S}$. Then for each  $(\eta,\vr,\tau)\in \overline {S}$, we can write
\be\label{presure-total-s}
h(\eta,\vr,\tau)  = h(\eta,\eta s_\vr,\eta s_\tau)  =   \de\,\eta^2 + k(L-1) \eta + a  \eta^\g s_\vr^\g  - \eta s_\tau, \quad s_\vr, s_\tau \in [0, \overline C].
\ee

 In $\overline {S}$, the monotonicity of the total pressure $h$ is mainly determined by $\eta$. We now show that even $h$ could be non monotone in $\eta$,  we can decompose it into a monotone part and a compactly supported part. Let $\overline R>1$ and $\chi\in C^\infty_c([0,\overline R))$ be a non-increasing cut-off function satisfying $0\leq \chi\leq 1$ and $\chi = 1$ on $ [0,\overline R_1]$, $0<\overline{R}_1 <\overline{R}$. We write the total pressure as
\be\label{H4-1}
 h(\eta,\eta s_\vr,\eta s_\tau)  = \mathcal{P}(\eta, s_\vr, s_\tau)  - \mathcal{R} (\eta,s_\vr, s_\tau),
\ee
with
\ba\label{H4-2}
 \mathcal{P}(\eta, s_\vr, s_\tau) & =  \de\,\eta^2 + k L\eta  + a  \eta^\g s_\vr^\g  - \big(1-\chi(\eta)\big) \big(  k \eta + \eta s_\tau  \big)   , \\
 \mathcal{R} (\eta,s_\vr, s_\tau)& =  \chi(\eta)\big(  k \eta + \eta s_\tau  \big) .
\ea
 By choosing $\overline R_1$ (and thus also $\overline{R}$) large enough, it is straightforward to check that the decomposition \eqref{H4-1}--\eqref{H4-2} satisfies Hypothesis (H4).

 \medskip

Hypothesis (H5) needs more attention, since it is closely connected to the form of the Helmholtz free energy. In fact, it is used in the construction of the approximate problem and it yields that $|\nabla^2_{\eta,\vr,\tau} H(\eta,\vr,\tau)| \leq C(r)(1+\eta^A)$ in the set $\{\eta^2 + \vr^2 +\tau^2>r^2\} \cap \overline{S}$. Hence, for our choice of the Helmholtz energy, we only need that
$$
|\nabla^2_\eta q(\eta)| + |\nabla^2_\vr p(\vr)| + |1/\tau|  \leq C(r)(1+\eta^A)
$$
in the set \{$\eta^2 + \vr^2 +\tau^2>r^2\} \cap \overline{S}$. However, it follows directly with the choice $A=0$ from the form of the pressure. The modified Hypothesis (H5) is fulfilled.

 \subsection{A priori estimates for smooth solutions}\label{sec:a priori}
In this section, we derive the a priori estimates for smooth solutions. Let $(\vr,\vu,\eta, \tau)$ be a smooth solution to \eqref{01}--\eqref{boundary} in $Q_T$ with smooth initial data satisfying the domination condition \eqref{ini-domination}. Moreover, without loss of generality, we assume the initial data for $(\vr,\eta, \tau)$ are bounded and strictly positive in $\overline\O$:
\be\label{ini-density-positive}
0<\underline a \leq \vr_0, \eta_0, \tau_0 \leq \overline a <\infty, \quad \mbox{for all $x\in \overline \O$}.
\ee
Otherwise, a standard trick is to mollify the initial data and add a strictly positive constant to the mollified data. Note that this can be done in such a way that the domination \eqref{ini-domination}  still holds.

 We start by showing the positivity of $\vr, \eta, \tau$ and the domination conditions \eqref{ini-domination-2}.
 \begin{alignat}{2}
\label{01-1}
\d_t \vr + \vu\cdot \nabla_x \vr &= (-\Div_x \vu)\vr,
\\
\label{03-1}
\d_t \eta + \vu\cdot \nabla_x \eta &= (-\Div_x \vu)\eta,
\\
\label{04-1}
\d_t  \tau + \vu\cdot \nabla_x \tau   &= \Big(-\Div_x \vu -  \frac{1}{2\lambda}\Big) \tau.
\end{alignat}
By introducing the characteristic along the velocity field
$$
\frac{\rm d}{\dt} X(t,x) = \vu (t,X(t,x)), \quad X(0,x) = x,
$$
the new unknowns along the characteristic defined by
$$
\tilde \vr(t,x) = \vr(t,X(t,x)), \quad \tilde \eta(t,x) = \eta(t,X(t,x)), \quad \tilde \tau(t,x) = \tau(t,X(t,x))
$$
solve
\begin{alignat}{2}
\label{01-2}
\d_t \tilde \vr  &= (-\Div_x \vu)(t,X(t,x))\tilde \vr,
\\
\label{03-2}
\d_t \tilde \eta  &= (-\Div_x \vu)(t,X(t,x))\tilde \eta,
\\
\label{04-2}
\d_t \tilde \tau  &= \Big(-\Div_x \vu -  \frac{1}{2\lambda}\Big)(t,X(t,x)) \tilde \tau.
\end{alignat}
It is straightforward to deduce that for all $(t,x)\in Q_T$,
\ba\label{density-positive}
 \inf_{x\in \O}\vr_0(x) \exp\left(-\int_0^t\|\Div_x \vu\|_{L^\infty(\O)}\,\dt'  \right) &\leq \vr(t,x) \leq \sup_{x\in \O}\vr_0(x) \exp\left(\int_0^t\|\Div_x \vu\|_{L^\infty(\O)}\,\dt'  \right), \\
 \inf_{x\in \O}\eta_0(x) \exp\left(-\int_0^t\|\Div_x \vu\|_{L^\infty(\O)}\,\dt'  \right) &\leq \eta(t,x) \leq \sup_{x\in \O}\eta_0(x) \exp\left(\int_0^t\|\Div_x \vu\|_{L^\infty(\O)}\,\dt'  \right), \\
 \inf_{x\in \O}\tau_0(x) \exp\left(-\int_0^t\|\Div_x \vu\|_{L^\infty(\O)}\,\dt' -  \frac{t}{2\lambda}  \right) &\leq \tau(t,x) \\
&\leq \sup_{x\in \O}\tau_0(x) \exp\left(\int_0^t\|\Div_x \vu\|_{L^\infty(\O)}\,\dt' -  \frac{t}{2\lambda} \right).
\ea
By \eqref{ini-density-positive} and \eqref{density-positive}, we know that $\vr,\eta, \tau$ are all bounded and strictly positive on $\overline Q_T$.

Similarly, for the differences $\xi: = \overline C \eta - \vr$ and $\zeta:  = \overline C \eta - \tau$, we have
\ba
&\d_t  \xi + \vu\cdot \nabla_x \xi = (-\Div_x \vu ) \xi, \quad &&\xi(0,\cdot) = \overline C \eta_0 - \vr_0 \geq 0, \\
&\d_t  \zeta + \vu\cdot \nabla_x \zeta   = (-\Div_x \vu )\zeta +  \frac{1}{2\lambda} \tau, \quad && \zeta(0,\cdot) = \overline C \eta_0 - \tau_0 \geq 0.
\ea
Direct calculation gives for all $(t,x)\in Q_T$:
\ba\label{density-positive-2}
  \xi(t,x) &\geq \inf_{x\in \O}\xi(0,x) \exp\left(\int_0^t - \|\Div_x \vu\|_{L^\infty(\O)}\,\dt'  \right) \geq 0, \\
 \zeta(t,x) &\geq \inf_{x\in \O}\zeta(0,x) \exp\left(\int_0^t - \|\Div_x \vu\|_{L^\infty(\O)}\,\dt'  \right) \\
& \quad +    \frac{t}{2\l} \inf_{x\in \O} \tau_0  \exp\left(\int_0^t -2 \|\Div_x \vu\|_{L^\infty(\O)}\,\dt' -  \frac{t}{2\lambda } \right)  \geq 0.
\ea

This implies that the domination condition preserves:
\be\label{ini-domination-3}
0\leq \vr(t,x) \leq \overline C \eta(t,x), \quad 0 \leq \tau(t,x) \leq \overline C \eta(t,x),\quad \mbox{for all $t,x \in \overline Q_T$}.
\ee
\medskip

Multiplying the momentum equation \eqref{02} by $\vu$ and using the other equations, one can derive the following energy equality
\ba\label{energy1}
&\int_\O \left[ \frac{1}{2} \vr |\vu|^2  + H(\vr,\eta, \tau) \right](t,\cdot)\dx  +  \int_0^t  \int_\O \SSS(\nabla_x \vu): \nabla_x \vu\,\dx \,\dt'  \\
& =  \int_\O \left[ \frac{1}{2} \vr_0 |\vu_0|^2 + H(\vr_0,\eta_0, \tau_0)  \right]\dx + \int_0^t \int_\O \vr\,\ff \cdot  \vu \,\dx\,\dt' + \frac{1}{2\lambda} \int_0^t \int_\O (\tau \log \tau + \tau)  \,\dx \,\dt',
\ea
where the Helmholtz free energy $H$ is given in \eqref{pressure-entropy}--\eqref{pressure-entropy-2}. As pointed in the introduction, the Helmholtz free energy may have a negative sign. While, this can be remedied by employing the domination condition \eqref{ini-domination-3}.  Indeed, if $\eta \geq \overline R_2$ for some $\overline R_2$ large, by using \eqref{ini-domination-3}, the quadratic term $\eta^2$ dominates and there holds $H(\vr, \eta, \tau) >0$. For $0\leq \eta \leq \overline R_2$, due to the continuity of $H$ and \eqref{ini-domination-3}, there exists a constant $\underline C>0$ such that $H(\vr, \eta, \tau) + \underline C >0. $ Hence, there exists a positive constant $\underline C$ such that
\be\label{pressure-entropy-3}
\tilde H (\vr, \eta, \tau) : = H(\vr, \eta, \tau) + \underline C >0,  \quad \mbox{for all $\eta \geq 0.$}
\ee
Again by \eqref{ini-domination-3}, there holds, by choosing $\underline C$ suitably large, that
$$
\tau^2 + \tau \log \tau + \tau \leq C \tilde H(\vr, \eta, \tau), \quad \vr +  \vr \log \vr  +  \vr^2 + \vr^\gamma \leq C \tilde H(\vr, \eta, \tau).
$$
We thus deduce from \eqref{energy1} that
\ba\label{energy2}
&\int_\O \left[ \frac{1}{2} \vr |\vu|^2  + \tilde H(\vr,\eta, \tau) \right](t,\cdot)\dx  +  \int_0^t  \int_\O \SSS(\nabla_x \vu): \nabla_x \vu\,\dx \,\dt'  \\
& \leq  \int_\O \left[ \frac{1}{2} \vr_0 |\vu_0|^2 + \tilde H(\vr_0,\eta_0, \tau_0)  \right]\dx + \frac{1}{2}\int_0^t \int_\O \vr|\vu|^2 \,\dx\,\dt' + \frac{C}{2\lambda} \int_0^t \int_\O \tilde H(\vr, \eta,\tau) \,\dx \,\dt'.
\ea
Hence, by Gronwall's  and Korn's inequalities, we obtain from \eqref{energy2} the following estimates
\ba\label{estimates-1}
\eta, \vr, \tau \in L^\infty (0,T; L^2(\Omega)), \quad \vu\in L^{2}(0,T;W_0^{1,2}(\Omega; \R^d)), \quad \vr |\vu|^2 \in L^\infty(0,T; L^{1}(\Omega)).
\ea

Next, in order to identify at least the weak limit of the pressure, we need to improve the $L^1$-estimate in the spatial variables. To this aim, the standard technique by Bogovskii-operator estimate can be employed, see \cite{bog} or the books \cite{Galdi-book, F-book, N-book}. Since this is a completely standard technique (see e.g. \cite[Formula (117)]{Novotny-Pokorny-19}), we do not present here any details. It yields, combined with the domination property,
$$
\eta, \vr, \tau \in L^{\frac 73} (Q_T).
$$

Finally, using \cite[Proposition 5]{Novotny-Pokorny-19} with a slight modification due to the additional damping term in the continuity equation for $\tau$ we easily verify that functions $s_\varrho$ and $s_\tau$, defined in \eqref{def-s} solve the following transport equations in the renormalized sense integrated up to the boundary and in the time integrated form (compare with \eqref{renormal-1})
\be\label{eq-s}
\d_t s_\vr + \vu \cdot \nabla_x s_\vr  = 0, \quad \d_t s_\tau + \vu \cdot \nabla_x s_\tau  =  - \frac{s_\tau}{2\l}.
\ee
Note that the assumption of Proposition 5 from \cite{Novotny-Pokorny-19} is fulfilled since all functions $\vr$, $\eta$ and $\tau$ are renormalized solutions to the continuity equation and therefore they are also in  $C([0,T];L^\th (\Omega))$ for all $1\leq \th <2$.

\subsection{Uniform estimates and basic convergence}

By using rather standard approach of constructing approximate solutions to compressible Navier-Stokes equations, that is to add artificial pressure terms, add diffusions for the continuity equations, and use Galerkin approximation for the momentum equations, one can get a family of approximate solutions which enjoy uniform estimates deduced from in Section \ref{sec:a priori}. We will not go through all the approximation schemes, because it can be done in the same manner as in Section 4 of \cite{Novotny-Pokorny-19}, and is similar as in \cite{F-book,FNP, N-book}.  Instead, we will start with $(\vr_n,\vu_n,\eta_n,\tau_n)_{n\in \N}$ which is a sequence of finite energy weak solutions in the sense of Definition \ref{def-weaksl-f} and we have the following uniform estimates deduced from the energy inequality

\ba\label{estimates-2}
&\sup_{n}\big(\|(\eta_n, \vr_n,\tau_n)\|_{L^\infty(0,T;L^2(\O))} + \|(\eta_n, \vr_n,\tau_n)\|_{L^{\frac 73}(Q_T)}\\
+ &\|\vu_n \|_{L^{2}(0,T;W_0^{1,2}(\Omega; \R^d))} + \|\vr_n |\vu_n|^2\|_{L^\infty(0,T; L^{1}(\Omega))} \big) < + \infty,
\ea
together with the domination relations
\be\label{ini-domination-4}
0\leq \vr_n(t,x) \leq \overline C \eta_n(t,x), \quad 0 \leq \tau_n(t,x) \leq \overline C \eta_n(t,x),\quad \mbox{for each $n\in \N$, for a.a. $(t,x) \in \overline Q_T$.}
\ee

We will pass to the limit in the sequel and show the compactness. The limit passage is done always up to a sequence, so we will not repeat this point.

By \eqref{estimates-2} using that $(\eta_n,\vr_n,\tau_n)$ solve the continuity equations, we have the weak convergence
\ba\label{limit-1}
&(\eta_n, \vr_n, \tau_n) \to (\eta, \vr, \tau) \ \mbox{weakly-* in} \ C_w([0,T];L^2(\O)) \quad \mbox{and weakly in} \ L^{\frac 73}(Q_T), \\
&  \vu_n\to \vu \ \mbox{weakly in} \ L^{2}(0,T;W_0^{1,2}(\Omega; \R^d)),
\ea
and the limit satisfies
\be\label{ini-domination-5}
0\leq \vr(t,x) \leq \overline C \eta(t,x), \quad 0 \leq \tau(t,x) \leq \overline C \eta(t,x),\quad \mbox{for a.a. $(t,x) \in \overline Q_T$.}
\ee
We use convention \eqref{def-s} to introduce
\be\label{rn-sn}
r_n = s_{\vr_n}  = \frac{\vr_n}{\eta_n} , \ s_n = s_{\tau_n} : = \frac{\tau_n}{\eta_n}.
\ee
Due to \cite[Proposition 5]{Novotny-Pokorny-19} we easily verify that functions $r_n$ and $s_n$ solve the following transport equations in the renormalized form integrated up to the boundary and in the time integrated form (see \eqref{renormal-3})
\be\label{eq-s-1}
\d_t r_n + \vu_n \cdot \nabla_x r_n  = 0, \quad \d_t s_n + \vu_n \cdot \nabla_x s_n  =  - \frac{s_n}{2\l}.
\ee
Due to the domination condition \eqref{ini-domination-5} on the limit, we can define by using the convention \eqref{def-s}
\be\label{def-r-s-tilde}
r : = \frac{\vr}{\eta}, \quad s : = \frac{\tau}{\eta}.
\ee
By \eqref{ini-domination-4} and \eqref{ini-domination-5}, there holds
\be\label{estimates-3}
0\leq r_n, s_n, r, s \leq \overline  C, \quad \forall\, n.\nn
\ee


Moreover, by using the continuity equations and the momentum equations, we have uniform estimates for the time derivatives $\d_t \vr_n, \d_t \eta_n, \d_t \tau_n$ and $\d_t(\vr_n \vu_n)$ in $L^\th(0,T;W^{-1,\th})$ for some $\th>1$. We thus have
\be\label{limit-3}
(\eta_n, \vr_n, \tau_n) \to (\eta, \vr, \tau) \ \mbox{in} \ C_{\rm weak}([0,T];L^2(\O)).
\ee
Together with $\vu_n\to \vu \ \mbox{weakly in} \ L^{2}(0,T;W_0^{1,2}(\Omega; \R^d))$ and by an Arzel\`a--Ascoli type argument we have the weak convergence for the nonlinear terms:
\be\label{limit-4}
\vr_n \vu_n \to \vr \vu, \  \eta_n \vu_n \to \eta\vu, \ \tau_n\vu_n \to \tau \vu, \ \vr_n \vu_n \otimes \vu_n \to \vr \vu \otimes \vu, \ \mbox{in $\mathcal{D}'(Q_T)$}.
\ee
Then the limit functions solve the continuity equations in the sense of distributions
 \begin{alignat}{2}
\label{01-3}
\d_t \vr + \Div_x (\vr \vu) &= 0,
\\
\label{03-3}
\d_t \eta + \Div_x (\eta \vu) &= 0,
\\
\label{04-3}
\d_t  \tau + \Div_x (\tau \vu) &=  -  \frac{1}{2\lambda} \tau.
\end{alignat}

\subsection{Renormalized continuity equations}\label{sec:renorm}

By applying Proposition 4 in \cite{Novotny-Pokorny-19} and using estimates \eqref{limit-1}, it is straightforward to deduce that all the continuity equations are satisfied in the renormalized sense up to the boundary:
 \ba\label{renormal-1}
\int_0^t \int_\O \big(b(\zeta_n) \d_t \phi + b(\zeta_n) \vu_n \cdot \nabla_x \phi - (b'(\zeta_n) \zeta_n - b(\zeta_n)) \Div_x \vu_n \, \phi\big) \,\dx\,\dt' \\
 = \int_\O b(\zeta_n) \phi (t,x) \,\dx  - \int_\O b(\zeta_n) \phi (0,x) \,\dx + \delta_{\zeta, \tau} \int_0^t \int_\O \frac{b'(\tau_n)\tau_n }{2\l}\,\dx\,\dt'
 \ea
and
 \ba\label{renormal-2}
\int_0^t \int_\O \big(b(\zeta) \d_t \phi + b(\zeta) \vu \cdot \nabla_x \phi - (b'(\zeta) \zeta - b(\zeta)) \Div_x \vu \, \phi\big) \,\dx\,\dt' \\
 = \int_\O b(\zeta) \phi (t,x) \,\dx  - \int_\O b(\zeta) \phi (0,x) \,\dx + \delta_{\zeta, \tau} \int_0^t \int_\O \frac{b'(\tau)\tau }{2\l}\,\dx\,\dt',
 \ea
for all $t\in [0,T]$ and for all $\phi \in C^1(\overline Q_T)$, where $\zeta$ and $\zeta_n$ can be arbitrarily chosen from $\{ \vr, \eta, \tau\}$ and $\{ \vr_n, \eta_n, \tau_n\}$, respectively. The equalities hold for all $b\in C([0,\infty)) \cap C^1((0,\infty))$ satisfying
\be\label{cond-b-1}
b'(\th) \th - b(\th) \in C[0,\infty), \  |b(\th)| + |b'(\th) \th - b(\th)| \leq C(1+\th) \ \mbox{for all $\th\geq 0.$}\nn
\ee
The constant $\delta_{\zeta, \tau}  = 1$ if $\zeta = \tau$; otherwise $\delta_{\zeta, \tau}  = 0$.

\medskip

For $r_n, s_n, r,s$, applying Proposition 5,  Proposition 7 (i.e., Proposition \ref{P2.5} in this paper), and Remark 3.2 in \cite{Novotny-Pokorny-19} implies that the transport equations \eqref{eq-s} hold in the renormalized sense up to the boundary:
 \ba\label{renormal-3}
&\int_0^t \int_\O \big(b(\zeta_n) \d_t \phi + b(\zeta_n) \vu_n \cdot \nabla_x \phi +  b(\zeta_n) \Div_x \vu_n \, \phi\big) \,\dx\,\dt' \\
&\quad = \int_\O b(\zeta_n) \phi (t,x) \,\dx  - \int_\O b(\zeta_n) \phi (0,x) \,\dx + \delta_{\zeta, r} \int_0^t \int_\O \frac{b'(r_n)r_n }{2\l}\,\dx\,\dt',
 \ea
and
 \ba\label{renormal-4}
&\int_0^t \int_\O \big(b(\zeta) \d_t \phi + b(\zeta) \vu \cdot \nabla_x \phi +  b(\zeta) \Div_x \vu \, \phi\big) \,\dx\,\dt' \\
&\quad = \int_\O b(\zeta) \phi (t,x) \,\dx  - \int_\O b(\zeta) \phi (0,x) \,\dx + \delta_{\zeta, r} \int_0^t \int_\O \frac{b'(r)r }{2\l}\,\dx\,\dt',
 \ea
for all $\phi \in C^1(\overline Q_T)$, where $\zeta$ and $\zeta_n$ can be arbitrarily chosen from $\{ r,s\}$ and $\{ r_n,s_n\}$, respectively. The equalities hold for any $b\in C[0,\infty) \cap C^1(0,\infty)$ satisfying
\be\label{cond-b-2}
b'(\th) s - b(\th) \in C[0,\infty), \  |b(\th)| \leq C(1+\th) \ \mbox{for all $\th\geq 0.$}\nn
\ee
The constant $\delta_{\zeta, r}  = 1$ if $\zeta = r$; otherwise $\delta_{\zeta, r}  = 0$.

\medskip


Applying again Proposition \ref{P2.5} and using its proof, we have for each $1\leq p <\infty$ that
\be\label{conv-r-s-tilde}
\int_\O \eta_n |r_n -  r|^p \, \dx \to 0, \quad \int_\O \eta_n |s_n -  s|^p \, \dx \to 0, \quad \mbox{for all $t\in [0, T]$}.
\ee
We remark that in \cite{VWY19}, a similar result is shown:
\be\label{conv-r-s-tilde-0}
\int_{0}^{T}\int_\O \eta_n |r_n -  r|^p \, \dx \to 0, \quad \int_{0}^{T}\int_\O \eta_n |s_n -  s|^p \, \dx \to 0.
\ee
The results in \eqref{conv-r-s-tilde} and \eqref{conv-r-s-tilde-0} offer some compactness for $r_{n}$ and $s_{n}$.
Since there is an extra damping term in the continuity equation in $\tau$, for the convenience of the readers, we briefly prove the second result in \eqref{conv-r-s-tilde}. 
 The strategy is to prove it for $p=2$ and the result for other $p$ follows from interpolation.  When $p=2$,
\be\label{conv-r-s-tilde-00}
\eta_n |s_n -  s|^2 = \eta_n s_n^2 + \eta_n  s^2 - 2 \eta_n s_n  s.
\ee
By  virtue of the renormalized equations in $s_n$ and $ s$, choosing $b(\th) = \th^{2}$ implies that  $s_n^2$ and $s^2$ satisfy
$$
\d_t s_n^2 + \vu \cdot \nabla_x s_n^2  =  - \frac{s_n^2}{\l}, \quad \d_t  s^2 + \vu \cdot \nabla_x s^2  =  - \frac{ s^2}{\l}.
$$
Proposition 6 in \cite{Novotny-Pokorny-19} implies that the products $\eta_n s_n^2$, $\eta_n s^2$ and $\eta  s^2$ satisfy the continuity equations  up to the boundary in the time integrated form:
$$
\d_t (\eta_n s_n^2) + \Div_x (\eta_n s_n^2 \vu)  =  - \frac{\eta_n s_n^2}{2\l}, \,\,  \d_t (\eta_n  s^2) + \Div_x (\eta_n  s^2 \vu)  =  - \frac{\eta_n   s^2}{2\l}, \,\, \d_t (\eta   s^2) + \Div_x (\eta   s^2 \vu)  =  - \frac{\eta    s^2}{2\l}.
$$
This gives
\ba\label{conv-r-s-tilde-1}
&\lim_{n\to \infty} \left(\int_\O (\eta_n s_n^2)(t)\,\dx  +  \int_0^t \int_\O  \frac{\eta_n s_n^2}{2\l} \,\dx\,\dt' \right)= \lim_{n\to \infty} \int_\O (\eta_n s_n^2)(0)\,\dx =  \int_\O \eta_0 s_0^2\,\dx,\\
&\lim_{n\to \infty} \left(\int_\O (\eta_n   s^2)(t)\,\dx  +  \int_0^t \int_\O  \frac{\eta_n   s^2}{2\l} \,\dx\,\dt' \right) = \lim_{n\to \infty} \int_\O (\eta_{n}    s^2)(0)\,\dx =  \int_\O \eta_0 s_0^2\,\dx,\\
&\int_\O (\eta   s^2)(t)\,\dx  +  \int_0^t \int_\O  \frac{\eta   s^2}{2\l} \,\dx\,\dt' =  \int_\O (\eta   s^2)(0)\,\dx =  \int_\O \eta_0 s_0^2\,\dx.
\ea
Again by Proposition 6 in \cite{Novotny-Pokorny-19}, the products $\tau_{n}   s$ and $\tau   s$ satisfy
$$
\d_t (\tau_n   s) + \Div_x (\tau_n   s \vu)  =  - \frac{\tau_{n}   s}{2\l}, \quad  \d_t (\tau   s) + \Div_x (\tau   s\vu)  =  - \frac{\tau    s}{2\l}.
$$
Then
\ba\label{conv-r-s-tilde-2}
& \lim_{n\to \infty} \left(\int_\O ( \tau_{n}   s)(t)\,\dx  +  \int_0^t \int_\O  \frac{\tau_{n}   s  }{2\l} \,\dx\,\dt' \right)  =  \lim_{n\to \infty} \int_\O ( \tau_{n}   s)(0)\,\dx = \int_\O \tau_0 s_0\,\dx,\\
&  \int_\O ( \tau   s)(t)\,\dx  +  \int_0^t \int_\O  \frac{\tau    s  }{2\l} \,\dx\,\dt'   =  \int_\O ( \tau   s)(0)\,\dx = \int_\O \tau_0 s_0\,\dx,
\ea
and furthermore
\ba\label{conv-r-s-tilde-3}
\lim_{n\to \infty} \left(\int_\O (\eta_n s_n   s)(t)\,\dx  +  \int_0^t \int_\O  \frac{\eta_n s_n   s }{2\l} \,\dx\,\dt' \right) &  = \lim_{n\to \infty} \left(\int_\O ( \tau_{n}   s)(t)\,\dx  +  \int_0^t \int_\O  \frac{\tau_{n}   s  }{2\l} \,\dx\,\dt' \right) \\
&  = \int_\O \tau_0 s_0\,\dx = \int_\O \eta_0 s_0^2\,\dx.
\ea
Hence, by \eqref{conv-r-s-tilde-00}--\eqref{conv-r-s-tilde-3}, we finally obtain
\be\label{conv-r-s-tilde-4}
\lim_{n\to \infty} \left(  \int_\O \eta_n |s_n -   s|^2(t) \, \dx + \int_0^t \int_\O  \frac{ \eta_n |s_n -   s|^2 }{2\l} \,\dx\,\dt' \right) = 0, \quad \mbox{for all $t\in [0,T]$}.
\ee
This implies \eqref{conv-r-s-tilde} by interpolation.

\subsection{Strong convergence of the densities}\label{sec:str-cov-density}

In order to show that
\be\label{w_p}
h(\eta_n,\vr_n,\tau_n) \to h(\eta,\vr,\tau) \quad \mbox{weakly in } L^1(Q_T),
\ee
we still need to show the strong convergence for $\eta_{n}$ and $\vr_{n}$. The proof is long and technical but  well understood nowadays: it can be done by employing the argument in Section 4.4 in \cite{Novotny-Pokorny-19}. So we only briefly recall the main steps, more details can be found in \cite{Novotny-Pokorny-19}.

\medskip

The starting point is to treat the total pressure as a solely function in $\eta_{n}$: we write $h(\eta_{n}, \vr_{n}, \tau_{n}) = h(\eta_{n}, \eta_{n} r_{n}, \eta_{n} s_{n})$ with $0\leq r_{n}, s_{n} \leq \overline C$. Then the idea is to use the nowadays well understood approach in the study of compressible Navier--Stokes equations (see for example \cite{FNP, F-book}) to prove the strong convergence of $\eta_{n}$.  Employing the effective viscous flux identity and using the renormalized equation in $\eta_{n}$ and $\eta$ implies that $\overline {\eta_{n} \log \eta_{n} } = \eta \log \eta$. The convexity of the function $\th \to \th \log \th$ in $[0,\infty)$ implies the strong convergence $\eta_{n} \to \eta$ strongly in $L^{1}(Q_{T})$. By interpolation, there holds $\eta_{n} \to \eta$ strongly in  $L^{\th}(Q_{T})$ for all  $1\leq \th<7/3$. This general idea, however, works only if the function  $\eta \mapsto h(\eta,\eta r,\eta \tau)$ is monotone which is not the case here. However, since it is monotone for large values of $\eta$ uniformly with respect to $r$ and $s$, we can, similarly as in  \cite[Sections 4.3 and 4.4]{Novotny-Pokorny-19}, apply the technique from \cite{Fe_non}.

\medskip

Let us now explain how we can get \eqref{w_p}, without having \eqref{eq8.7} from Hypothesis (H3).
We first apply \eqref{conv-r-s-tilde}-\eqref{conv-r-s-tilde-0} with $p=1$ to deduce
$$
\eta_{n} r_{n} - \eta_{n} r \to 0, \  \eta_{n} s_{n} - \eta_{n} s \to 0, \  \mbox{in}  \ L^{1}(Q_{T}).
$$
Together with the strong convergence $\eta_{n} \to \eta$ in $L^{1}(Q_{T})$,  we finally obtain
\ba\label{strong-conv}
& \vr_{n} - \vr = \eta_{n} r_{n}  - \eta r = (\eta_{n} r_{n}  - \eta_{n} r) + (\eta_{n} - \eta) r \to 0, \\
& \tau_{n} - \tau = \eta_{n} s_{n}  - \eta s = (\eta_{n} s_{n}  - \eta_{n} s) + (\eta_{n} - \eta) s \to 0,\nn
\ea
 strongly in $L^{1}(Q_{T})$, hence (up to a subsequence) a.e. in $Q_T$. Furthermore, the convergence is also strong in $L^{\th}(Q_{T})$ for each $1<\th<\frac 73$ by interpolation. Thus, by virtue of the Vitali's convergence theorem we conclude \eqref{w_p}.  This allows us to pass to the limit  and deduce the limit momentum equations.

 \subsection{Energy inequality and end of the proof}
To finish the proof, it remains to show the energy inequality.  We recall the energy inequality for $(\vr_{n}, \vu_{n}, \eta_{n}, \tau_{n})$: for a.a. $t \in (0,T)$,
\ba\label{energy1-f-n}
&\int_\O \left[ \frac{1}{2} \vr_{n} |\vu_{n}|^2  + H(\vr_{n},\eta_{n}, \tau_{n}) \right](t,\cdot)\dx  +  \int_0^t  \int_\O \SSS(\nabla_x \vu_{n}): \nabla_x \vu_{n}\,\dx \,\dt'  \\
& \leq  \int_\O \left[ \frac{1}{2} \vr_0 |\vu_0|^2 + H(\vr_0,\eta_0, \tau_0)  \right]\dx + \int_0^t \int_\O \vr\,\ff \cdot  \vu_{n} \,\dx\,\dt' + \frac{1}{2\lambda} \int_0^t \int_\O (\tau_{n} \log \tau_{n} + \tau_{n})  \,\dx \,\dt'.
\ea
The Helmholtz free energy $H$ has the form \eqref{pressure-entropy}--\eqref{pressure-entropy-2}. The idea is to pass $n\to \infty$ in \eqref{energy1-f-n} to get the energy inequality for the limit $\vr, \vu, \eta, \tau.$ Compared to the compressible Navier--Stokes equations, we need to deal with the extra integrals:
\ba\label{energy-ineq-1}
k(L-1) \int_\O \eta_n \log \eta_n \,\dx , \quad  - \int_\O  \tau_n \log\tau_n \,\dx , \quad \int_\O  \vr_n^\g\,\dx , \ 0<\g<1.
\ea
The other terms can be dealt as in the study for the compressible Navier--Stokes equations, by using the weak lower semicontinuity and the convergences in \eqref{limit-1} and in Section \ref{sec:str-cov-density}, and the fact that $\th^2, \ \th^\g$ with $\g\geq 1$ and  $\th \log \th$ are all convex functions in $\th$ on $(0,\infty)$.

The issues for the terms in \eqref{energy-ineq-1} are respectively that the coefficient $L-1$ is possibly negative and $-1$ is negative, and the function $\vr \to \vr^\g$ is not convex when $0<\g<1.$  The idea is to consider the renormalized equations:
\ba\label{weak-renormal-123}
& \d_t (\vr_n^\g) +\Div_x (\vr_n^\g\vu_n) + (\g-1)\vr_n^\g\,\Div_x \vu_n =0,\\
& \d_t (\eta_n \log \eta_n) +\Div_x (\eta_n\log \eta_n \vu) + \eta_n \,\Div_x \vu_n =0, \\
& \d_t (\tau_n \log \tau_n) +\Div_x (\tau_n\log \tau_n \vu) + \tau_n \,\Div_x \vu_n = -\frac{ \tau_n(\log\tau_n+1)}{2\l},\nn
\ea
which hold by using the estimates in \eqref{estimates-2} and observing that
$$\vr_n^\g\,\Div_x \vu_n, \ \  \eta_n \,\Div_x \vu_n, \ \  \tau_n \,\Div_x \vu_n  + \frac{ \tau_n(\log\tau_n+1)}{2\l} $$
are weakly convergent to their counterparts (removing the lower index $n$) in $L^1(Q_T)$.  Then again by the estimates in \eqref{estimates-2},  we have
$$
\vr_n^\g, \  \eta_n\log\eta_n,\ \tau_n\log\tau_n \in C_w([0,T],L^1(\O)).
$$
Together with the strong convergence we have shown in Section \ref{sec:str-cov-density} and the estimates from Section 3.4, we finally have
$$
\vr_n^\g \to \vr^\g, \  \eta_n\log\eta_n \to \eta \log \eta,\ \tau_n\log\tau_n \to \tau \log\tau  \ \mbox{in} \  C_w([0,T],L^1(\O)).
$$
Then passing to the limit $n\to \infty$ gives our desired energy inequality for the limit solution $(\vr, \vu, \eta, \tau)$ and we complete the proof.

\medskip

{\bf Acknowledgement:} The work of Y.L. has been supported by the Recruitment Program of Global Experts of China. The work of M.P. was supported by the grant of the Czech Science Foundation No. 19-04243S. The paper was written during the stay of the second author at Nanjing University. This paper is dedicated to Prof. Weiyi Su on her 80 years birthday.


\end{document}